# An unconstrained-like control-based dynamic method for optimization problems with simple bounds

Sheng ZHANG, Xin DU, Fang-Fang HU, and Jiang-Tao HUANG

**Abstract**: The optimization problems with simple bounds are an important class of problems. To facilitate the computation of such problems, an unconstrained-like dynamic method, motivated by the Lyapunov control principle, is proposed. This method employs the anti-windup limited integrator to address the bounds of parameters upon the dynamics for unconstrained problem, and then solves the transformed Initial-value Problems (IVPs) with mature Ordinary Differential Equation (ODE) integration methods. It is proved that when the gain matrix is diagonal, the result is equivalent to that of the general dynamic method which involves an intermediate Quadratic Programming (QP) sub-problem. Thus, the global convergence to the optimal solution is guaranteed without the requirement of the strict complementarity condition. Since the estimation of the right active constraints is avoided and no programming sub-problem is involved in the computation process, it shows higher efficiency than the general dynamic method and other common iterative methods through the numerical examples. In particular, the implementation is simple, and the proposed method is easy-to-use.

## 1. Introduction

The optimization problems with simple bounds, also called the box constrained optimization problems, is a class of optimization problems with a special structure. The mathematical model of these problems is as follows:

***Problem 1***: Consider the objective function

$$J = f(\boldsymbol{\theta}), \qquad (1)$$

subject to

$$\boldsymbol{\theta}^l \leq \boldsymbol{\theta} \leq \boldsymbol{\theta}^u, \qquad (2)$$

This paper is based upon work supported by the National Natural Science Foundation under Grant No. 11902332. The authors are with the China Aerodynamics Research and Development Center, Mianyang, 621000, China. (email: zszhangshengzs@hotmail.com).

where $\boldsymbol{\theta} \in \mathbb{R}^n$ is the optimization parameter vector. $f : \mathbb{R}^n \to \mathbb{R}$ is a scalar function with continuous first-order partial derivatives with respect to $\boldsymbol{\theta}$. $\boldsymbol{\theta}^l$ and $\boldsymbol{\theta}^u$ are the lower bounds and the upper bounds of the parameters, respectively. Find the optimal solution $\hat{\boldsymbol{\theta}}$ that minimizes $J$. Note that for the unconstrained parameter $\theta_i$, it may be set that $(\boldsymbol{\theta}^l)_i = -\inf$ and $(\boldsymbol{\theta}^u)_i = \inf$.

The simple bound constrained problems may be found in many practical applications. Actually, it is said that any real-world unconstrained problem is meaningful only when considered with the bound constraints [1]. Meanwhile, Such problem is often a sub-problem of the augmented Lagrangian or of the penalty computational schemes for solving general constrained optimization [2, 3]. Therefore, the development of numerical algorithms to efficiently solve Problem 1 is important in both theory and practice [4].

The approaches most widely used to solve Problem 1 are the active-set methods [5-8]. In such method, at each iteration we have a working set that is supposed to approximate the set of active constraints and the set is iteratively updated. The main difficulty with these method is the determination of the right active sets in an iterative process either by adding parameters which violate one of their bounds or by removing those for which further progress is predicted [4]. In general, only a single active constraints can be added or deleted to the active set at each iteration unless the strict complementarity condition holds [6], and this makes the computation slow. There is a better understanding towards such method and efficient algorithms have been developed with desirable performance [9-12].

The trust-region type algorithms are also employed to address the box constrained optimization problems [13-16]. Using such method, the optimal solution may be obtained globally [15], and a superlinear convergence rate may be established without the strict complementarity condition under appropriate assumptions [16]. Upon the trust-region method, the interior-point methods for Problem 1 are further developed and they are proven to be robust and efficient [17, 18]. Interesting theoretical results are obtained in this framework and superlinear or even quadratic convergence is obtained without the strict complementarity assumption [19]. Because Problem 1 falls into the sub-category of the general constrained Non-linear Programing (NLP) problems, it may be solved by the Sequential Quadratic Programming (SQP) method, which is among the most prevailing and efficient methods for the constrained optimization [1].

Most of the methods mentioned above involve an immediate problem, either for programing a sub-problem or for determining a proper step size. To enhance the efficiency, many researches exploit the unconstrained optimization upon the special structure of Problem 1 when devising the algorithms. The gradient projection method takes advantage of the gradient processed by the projection operator, and may be almost as simple as their unconstrained counterparts [6, 20]. The unconstrained optimization is combined to effectively exploit the box structure of the constraints during the computation [3, 7]. With a differentiable penalty function, an unconstrained minimization approach for the box constrained problem is proposed [21]. Besides, Problem 1 may be directly reformed as an unconstrained problem to solve, by transforming the parameters. Concretely, quadratic transformation for parameters with single bounds and "sin" transformation for parameters with dual bounds [22].

Besides the iterative methods that use discrete mechanism, there is another way to solve the optimization problems, which is based on the continuous-time dynamics. With such method, the Dynamic Optimization Equation (DOE), which is a set of differential equations, is developed and the optimization problem is transformed to the Initial-value Problem (IVP) to be solved. Studies on this type of methods may date back to the 1940s for the unconstrained problems [23]. Relevant work on the unconstrained minimization includes the gradient dynamic equation [24-26], the second-order dynamic equation arising from the physical energy view [27], and the continuous Newton method [28], etc. To address the constrained NLP problems, the DOE for the NLP with ECs are developed and the quadratic slack parameters are employed to address the inequality constraints [29-32]. The DOE, defined on the original optimization parameters, for the NLP with IECs is also proposed in the light of Riemannian metric [33]. Under the framework of the interior-point method, the DOE proposed may also avoid the introduction of slack parameters; see [34] and [35]. Recently, a control-based dynamic method for the general NLPs, based on the Lyapunov stability principle, is proposed [36].The DOE, which has the same dimension to that of the optimization parameter vector, is established and its solution will converge to the optimal solution of the NLP globally, regardless of the feasibility of the initial condition.

In this paper we propose an unconstrained-like dynamic method to solve Problem 1. This method also arises from the control-based view but solves Problem 1 with an unconstrained-like strategy. It simply uses the limited integrator to address the bounds, and therefore is easy-to-use. The anti-windup limited integrator acting on the continuous dynamics has a similar role to the projection operator in the iterative gradient projection method. However, the daunting task of searching for a reasonable step size and the annoying oscillation phenomenon around the optimum are eliminated. The paper is organized as follows. First the optimality condition for the box constrained problems is presented in Section 2. Then Section 3 gives the general control-based dynamic method and its concrete form when applying to the optimization problems with simple bounds. In Section 4, the unconstrained-like dynamic method is established, including the theoretical foundation and the numerical issues. Later in Section 5, numerical examples are solved to verify the proposed method. Section 7 concludes the article at last.

## 2. Optimality conditions of Problem 1

First necessary definitions on the NLP are given for the following study. The general NLP may be formulated as

$$\begin{aligned} \min \quad & J = f(\boldsymbol{\theta}) \\ \text{subject to} \quad & \\ & \boldsymbol{g}(\boldsymbol{\theta}) \leq \boldsymbol{0} \\ & \boldsymbol{h}(\boldsymbol{\theta}) = \boldsymbol{0} \end{aligned} \quad , \tag{3}$$

where $\boldsymbol{\theta} \in \mathbb{R}^n$ is the optimization parameter vector. $f : \mathbb{R}^n \to \mathbb{R}$ is a scalar function with continuous first-order partial derivatives with respect to $\boldsymbol{\theta}$. $\boldsymbol{g} : \mathbb{R}^n \to \mathbb{R}^r$ and $\boldsymbol{h} : \mathbb{R}^n \to \mathbb{R}^s$ are $r$-dimensional vector function and $s$-dimensional vector function with continuous first-order partial derivatives, respectively.

**Definition 1:** The feasible solution region $\mathbb{D}_f$ for the NLP is the collection of parameters that satisfy the equality and inequality constraints, i.e., $\mathbb{D}_f = \{\boldsymbol{\theta} \mid \boldsymbol{g}(\boldsymbol{\theta}) \leq \boldsymbol{0}, \boldsymbol{h}(\boldsymbol{\theta}) = \boldsymbol{0}\}$.

With Definition 1, the NLP defined in Eq. (3) may be expressed concisely, namely

$$\hat{\boldsymbol{\theta}} = \underset{\boldsymbol{\theta} \in \mathbb{D}_f}{\arg \min}(J) . \tag{4}$$

**Definition 2:** For the $i$-th inequality constraint in Eq. (3), it is said to be an optimal active inequality constraint if $g_i(\hat{\boldsymbol{\theta}}) = 0$.

Note that the optimal active inequality constraint is defined on the optimal solution $\hat{\boldsymbol{\theta}}$, while other inequality constraints may still be activated or even violated during the optimization process. Problem 1 is a special case of the general constrained problem. In order to facilitate the study, Problem 1 may be rearranged to the standard form as follows

$$\begin{aligned} \min \quad & J = f(\boldsymbol{\theta}) \\ \text{subject to} \quad & \\ & \boldsymbol{g}(\boldsymbol{\theta}) \leq \boldsymbol{0} \end{aligned} \quad , \tag{5}$$

Here concretely $\boldsymbol{g}(\boldsymbol{\theta}) = \begin{bmatrix} \boldsymbol{\theta} - \boldsymbol{\theta}^u \\ -\boldsymbol{\theta} + \boldsymbol{\theta}^l \end{bmatrix}$. It is well know that the Karush-Kuhn-Tucker (KKT) optimality condition characterizes the optimal solution and provides the basis for algorithms. In practice, its existence is ensured by the Constraint Qualification (CQ) in that the linear approximation of the constraints at the optimal solution adequately captures the geometry of the feasible set. For Problem 1, the constraint qualification holds since the constraints are all linear. Therefore, the KKT optimality conditions always exists [2]; that is

$$\boldsymbol{f}_{\boldsymbol{\theta}} + \boldsymbol{g}_{\boldsymbol{\theta}}^{\mathrm{T}} \boldsymbol{\pi}_I = \boldsymbol{0} , \tag{6}$$

$$\begin{aligned} (\pi_I)_i \geq 0 & \quad i \in \hat{\mathbb{I}} \\ (\pi_I)_i = 0 & \quad i \notin \hat{\mathbb{I}} \end{aligned} , \tag{7}$$

where $f_{\boldsymbol{\theta}} = \dfrac{\partial f}{\partial \boldsymbol{\theta}}$ is the shorthand notation in the form of column vector. $\boldsymbol{g}_{\boldsymbol{\theta}}$ are the Jacobi matrixes. $\hat{\mathbb{I}}$ is the index set of the optimal active inequality constraints for Problem 1, defined as

$$\hat{\mathbb{I}} = \{i \mid g_i(\hat{\boldsymbol{\theta}}) = 0, \ i = 1, 2, ..., r\}. \tag{8}$$

Eq. (7) may be re-presented as

$$(\pi_I)_i g_i = 0, \tag{9}$$

and it is called the complementarity condition. The strict complementarity condition requires that

$$(\pi_I)_i > 0 \quad i \in \hat{\mathbb{I}}. \tag{10}$$

Obviously, when the optimal parameter $\hat{\theta}_i$ does not reach the bound, we have $f_{\hat{\theta}_i} = 0$. In particular, the Jacobi matrix of $\boldsymbol{g}_{\boldsymbol{\theta}}$ is $\boldsymbol{g}_{\boldsymbol{\theta}} = \begin{bmatrix} I_{n \times n} \\ -I_{n \times n} \end{bmatrix}_{2n \times n}$, then for the parameter which is active at the upper bound, we have $(\pi_I)_i = -(f_{\boldsymbol{\theta}})_i$, and if it is active at the lower bound, we have $(\pi_I)_i = (f_{\boldsymbol{\theta}})_i$. If the optimal parameter $\hat{\theta}_i$ take the value of the bound and $f_{\hat{\theta}_i} = 0$ meanwhile, the corresponding multiplier is zero. Then the strict complementarity condition is violated and the corresponding solution $\hat{\boldsymbol{\theta}}$ is called degenerate.

## 3. The control-based dynamic method

### 3.1. The control-based view upon Lyapunov stability theory

The Lyapunov stability theory in the control field investigates the dynamic behaviour of states within a dynamic system, from the view of generalized energy [37]. It states that for a stable dynamic system

$$\dot{\boldsymbol{x}} = \boldsymbol{f}(\boldsymbol{x}), \tag{11}$$

where $\boldsymbol{x}$ is the state, $\dot{\boldsymbol{x}} = \frac{d\boldsymbol{x}}{dt}$ is its time derivative, and $\boldsymbol{f} : \mathbb{D} \to \mathbb{R}^n$ is a vector function, a monotonically decreasing 'energy' function $V(\boldsymbol{x})$, which achieves its minimum when the state $\boldsymbol{x}$ reaches the equilibrium point $\hat{\boldsymbol{x}}$, may be constructed. The Lyaounov principle is also the foundation for the system synthesis. Based on the control-based view, a dynamic method to solve the general NLP problem is proposed in Ref. [36], through seeking the dynamics that minimize the objective function (regarded as the measure of generalized energy).

Consider the dynamics

$$\frac{d\boldsymbol{\theta}}{d\tau} = \boldsymbol{u}, \tag{12}$$

where $\tau$ is the virtual time introduced to describe the expected dynamics and $\boldsymbol{u}$ is the control to be determined. The motivation is to find $\boldsymbol{u}$ that makes the optimal solution the stable equilibrium point of the resulting dynamic system. Then the optimal solution may be solved in an asymptotic manner. This problem is a bit like the controller design. However, here the equilibrium is unknown while the Lyapunov function already exists. Following such idea, the DOE, which has the same dimension to that of the optimization parameter vector, may be established and its solution will converge to the optimal solution of the NLP globally. Via this control-based method, the NLP may be transformed to the IVP to be solved, with mature Ordinary Differential Equation (ODE) integration methods. In contrast to the iterative methods, in which the determination of a proper step size is an artful task, the dynamic method avoids this difficulty and the possible numerical oscillation around the optimum, by employing the mature ODE integration methods. Various reliable numerical approaches, including the explicit Runge-Kutta methods, the implicit Back Differentiation Formulas (BDF) methods and other variable-order variable-step techniques, etc. [38], may be applied to solve the resulting IVP conveniently.

### 3.2. The unconstrained dynamic method

With the dynamic method, an unconstrained problem may be easily solved as follows. To find its optimal value $\hat{\boldsymbol{\theta}}$ that minimizes the performance index $J$, an analogy with the Lyapunov function is made and $J$ is differentiated with respect to $\tau$ to give

$$\frac{\mathrm{d}J}{\mathrm{d}\tau} = \frac{\mathrm{d}f}{\mathrm{d}\tau} = f_{\boldsymbol{\theta}}^{\mathrm{T}} \boldsymbol{u} . \tag{13}$$

To guarantee that $J$ decreases with respect to $\tau$, i.e., $\frac{\mathrm{d}J}{\mathrm{d}\tau} \leq 0$, the control law may be set as

$$\boldsymbol{u} = -\boldsymbol{K}_{\boldsymbol{\theta}} f_{\boldsymbol{\theta}} , \tag{14}$$

and the resulting DOE is

$$\frac{\mathrm{d}\boldsymbol{\theta}}{\mathrm{d}\tau} = -\boldsymbol{K}_{\boldsymbol{\theta}} f_{\boldsymbol{\theta}} , \tag{15}$$

where $\boldsymbol{K}_{\boldsymbol{\theta}}$ is an $n \times n$ dimensional positive-definite gain matrix. In seeking the optimal solution, the initial condition for Eq. (15) is allowed to be arbitrary. Note that here $\boldsymbol{K}_{\boldsymbol{\theta}}$ is only required to be positive-definite and it may be either constant or dependent on $\boldsymbol{\theta}$. A constant case is simple to use, while a parameter-dependent one may bring better performance. However, the tuning will be an artful task.

### 3.3. The general dynamic method for Problem 1

For Problem 1, determining a feasible solution is easy. Therefore, the dynamic method operating in the feasible solution region will be employed. In applying the standard procedure in [36] for the control law $\boldsymbol{u}$, first Problem 1 is rearranged to the standard form, as given in Eq. (5). Note that the differential constraints arising from the IECs $\boldsymbol{g}(\boldsymbol{\theta}) \leq \boldsymbol{0}$ need to be considered in contrast to the unconstrained problems. Then the Feasibility-preserving Dynamics Optimization Problem (FPDOP), which is a Quadratic Programming (QP) problem within the feasible solution region, was constructed to seek $\boldsymbol{u}$; that is

$$\begin{aligned} \min \quad & J_{t3} = \frac{1}{2} J_{t1} + \frac{1}{2} J_{t2} \\ \text{subject to} \quad & \\ & (g_i)_{\boldsymbol{\theta}}^{\mathrm{T}} \boldsymbol{u} \leq 0 \quad i \in \mathbb{I} \end{aligned} , \tag{16}$$

where $\mathbb{I}$ is the index set of the activated IECs for the parameters $\boldsymbol{\theta}$ and it is defined as

$$\mathbb{I} = \{i \mid g_i(\boldsymbol{\theta}) = 0, \ i = 1, 2, ..., 2n\}. \tag{17}$$

The performance index $J_{t1}$ and $J_{t2}$ are used to reduce the objective function and avoid the infinite control, respectively. They are given by

$$J_{t1} = f_{\boldsymbol{\theta}}^{\mathrm{T}} \boldsymbol{u}, \tag{18}$$

$$J_{t2} = \frac{1}{2} \boldsymbol{u}^{\mathrm{T}} \boldsymbol{K}_{\boldsymbol{\theta}}^{-1} \boldsymbol{u}. \tag{19}$$

Introduce the index set of optimal active IECs for the FPDOP, namely

$$\mathbb{I}_p = \{i \mid g_i = 0, \ (g_i)_{\boldsymbol{\theta}}^{\mathrm{T}} \boldsymbol{u} \leq 0 \text{ is an optimal active IEC for the FPDOP (16)}, \ i = 1, 2, ..., r\}, \tag{20}$$

The number of elements in $\mathbb{I}_p$ is denoted by $n_{\mathbb{I}_p}$. In particular, for the optimal solution $\hat{\boldsymbol{\theta}}$, $\mathbb{I}_p$ is highlighted with a hat '^' as $\hat{\mathbb{I}}_p$. Now the following variant of the FPDOP (16) is presented as

$$\begin{aligned} \min \quad & J_{t3} = \frac{1}{2} J_{t1} + \frac{1}{2} J_{t2} \\ \text{subject to} \quad & \\ & (g_i)_{\boldsymbol{\theta}}^{\mathrm{T}} \boldsymbol{u} = 0 \quad i \in \mathbb{I}_p \end{aligned} \tag{21}$$

To derive this form, the key issue is to determine the index set $\mathbb{I}_p$ from $\mathbb{I}$, which may be obtained by solving the QP sub-problem (16). Once ascertained, the DOE for Problem 1 may be derived. By applying the general results in Ref. [36] to Problem 1, Theorem 1 below may be established.

**Lemma 1** [36]: *For the general NLP given by Eq. (3), suppose the CQ holds so that the KKT optimality condition exists. Solve the IVP defined by the following DOE*

$$\frac{\mathrm{d}\boldsymbol{\theta}}{\mathrm{d}\tau} = -\boldsymbol{K}_{\boldsymbol{\theta}}(f_{\boldsymbol{\theta}} + \boldsymbol{h}_{\boldsymbol{\theta}}^{\mathrm{T}} \boldsymbol{\pi}_E + \boldsymbol{g}_{\boldsymbol{\theta}}^{\mathrm{T}} \boldsymbol{\pi}_I), \tag{22}$$

*with arbitrary initial condition* $\boldsymbol{\theta}|_{\tau=0} \in \mathbb{D}_f$, *where the parameter vectors* $\boldsymbol{\pi}_E \in \mathbb{R}^s$ *and* $\boldsymbol{\pi}_I \in \mathbb{R}^r$ *are determined by*

$$\boldsymbol{\pi}_E = \begin{bmatrix} \pi_1 \\ \pi_2 \\ \vdots \\ \pi_s \end{bmatrix}, \quad (\pi_I)_i = 0 \quad i \notin \mathbb{I}_p, \quad \boldsymbol{\pi}_I(\mathbb{I}_p) = \begin{bmatrix} \pi_{s+1} \\ \pi_{s+2} \\ \vdots \\ \pi_{s+n_{\mathbb{I}_p}} \end{bmatrix} \geq \boldsymbol{0} \tag{23}$$

*and the parameter vector* $\boldsymbol{\pi} \in \mathbb{R}^{s+n_{\mathbb{I}_p}}$ *is computed by*

$$\boldsymbol{\pi} = -(\tilde{\boldsymbol{h}}_\theta \boldsymbol{K}_\theta \tilde{\boldsymbol{h}}_\theta^{\mathrm{T}})^+ \tilde{\boldsymbol{h}}_\theta \boldsymbol{K}_\theta \boldsymbol{f}_\theta. \tag{24}$$

$\boldsymbol{K}_\theta$ is an $n \times n$ dimensional positive-definite gain matrix and $\tilde{\boldsymbol{h}} = \begin{bmatrix} \boldsymbol{h} \\ \boldsymbol{g}(\mathbb{I}_p) \end{bmatrix}$. $\mathbb{I}_p$ is the index set of optimal active inequality constraints for the FPDOP. When $\tau \to +\infty$, $\boldsymbol{\theta}$ will satisfy the KKT optimality condition. Moreover, for the optimal solution, $\hat{\mathbb{I}}_p = \hat{\mathbb{I}}$ and the value of $\boldsymbol{\pi}$ is independent of $\boldsymbol{K}_\theta$.

**Theorem 1:** *For Problem* 1*, solve the IVP defined by the following DOE*

$$\frac{\mathrm{d}\boldsymbol{\theta}}{\mathrm{d}\tau} = -\boldsymbol{K}_\theta (\boldsymbol{f}_\theta + \boldsymbol{g}_\theta^{\mathrm{T}} \boldsymbol{\pi}_I), \tag{25}$$

*with arbitrary initial condition* $\boldsymbol{\theta}|_{\tau=0} \in \mathbb{D}_f$*, where the parameter vectors* $\boldsymbol{\pi}_I \in \mathbb{R}^{2n}$ *is determined by*

$$(\boldsymbol{\pi}_I)_i = 0 \quad i \notin \mathbb{I}_p, \quad \boldsymbol{\pi}_I(\mathbb{I}_p) = \begin{bmatrix} \pi_1 \\ \pi_2 \\ \dots \\ \pi_{n_{\mathbb{I}_p}} \end{bmatrix} \geq \boldsymbol{0}, \tag{26}$$

*and the parameter vector* $\boldsymbol{\pi} \in \mathbb{R}^{n_{\mathbb{I}_p}}$ *is computed by*

$$\boldsymbol{\pi} = -(\tilde{\boldsymbol{h}}_\theta \boldsymbol{K}_\theta \tilde{\boldsymbol{h}}_\theta^{\mathrm{T}})^{-1} \tilde{\boldsymbol{h}}_\theta \boldsymbol{K}_\theta \boldsymbol{f}_\theta. \tag{27}$$

$\boldsymbol{K}_\theta$ *is an* $n \times n$ *dimensional positive-definite gain matrix and* $\tilde{\boldsymbol{h}} = \boldsymbol{g}(\mathbb{I}_p)$. $\mathbb{I}_p$ *is the index set of optimal active inequality constraints for the FPDOP. When* $\tau \to +\infty$*,* $\boldsymbol{\theta}$ *will satisfy the optimality conditions* (6) *and* (7)*. Moreover, for the optimal solution,* $\hat{\mathbb{I}}_p = \hat{\mathbb{I}}$ *and the value of* $\boldsymbol{\pi}$ *is independent of* $\boldsymbol{K}_\theta$*.*

In pursuing the optimal solution under the dynamics governed by Eq. (25), an intermediate QP sub-problem will be solved and the set $\mathbb{I}_p$ needs to be determined dynamically. Generally, $\mathbb{I}$ is easy to get and is a good initial guess for $\mathbb{I}_p$. Thus all inequality constraints in $\mathbb{I}$ may be first strengthened to get the corresponding Lagrange multipliers and then use the optimality conditions for the QP to seek the right $\mathbb{I}_p$, following the procedure of the well-known active-set methods [1, 2].

In particular, the $n_{\mathbb{I}_p} \times n$ dimensional Jacobi matrix $\tilde{\boldsymbol{h}}_\theta$ has full row rank and each row only has one nonzero element of 1 or -1. It is interesting to find that $(\tilde{\boldsymbol{h}}_\theta)^+ = (\tilde{\boldsymbol{h}}_\theta)^{\mathrm{T}}$ and $\tilde{\boldsymbol{h}}_\theta \tilde{\boldsymbol{h}}_\theta^+ = \tilde{\boldsymbol{h}}_\theta \tilde{\boldsymbol{h}}_\theta^{\mathrm{T}} = \boldsymbol{I}$, where the superscript "+" denotes the pseudo-inverse of matrix and $\boldsymbol{I}$ is the identity matrix. In contrast, $\tilde{\boldsymbol{h}}_\theta^+ \tilde{\boldsymbol{h}}_\theta$ is a square matrix, where only the diagonal elements may take the value of 1 if the corresponding column vector in $\tilde{\boldsymbol{h}}_\theta$ is not a zero vector.

For example, if $\tilde{\boldsymbol{h}}_\theta = \begin{bmatrix} 1 & 0 & 0 & 0 \\ 0 & 0 & -1 & 0 \end{bmatrix}$, then $\tilde{\boldsymbol{h}}_\theta^+ \tilde{\boldsymbol{h}}_\theta = \begin{bmatrix} 1 & 0 & 0 & 0 \\ 0 & 0 & 0 & 0 \\ 0 & 0 & 1 & 0 \\ 0 & 0 & 0 & 0 \end{bmatrix}$. Likewise, Theorem 1 implies the convergence of the multipliers, namely

**Remark 1:** Solve Problem 1 with the approach suggested in Theorem 1; the KKT multiplier given by Eqs. (26) and (27) will converge to the right solution as $\theta$ converges to $\hat{\theta}$.

## 4. The proposed unconstrained-like dynamic approach

### 4.1. Theoretical foundation

The dynamic method in the last section may solve the general NLP. However, an intermediate FPDOP is involved, which increases the computation burden. For Problem 1, the method proposed in the following is more easy-to-use. It treats the problem like an unconstrained problem and directly use the limited integration to address the parameter bounds, which is a common technique for anti-windup in the integral control. Thus, no extra programming problem needs to be solved during the computation progress, and this may improve the computation efficiency.

As given in Sec. 3.2, an unconstrained problem may be solved with the dynamic method, and the integral form for Eq. (15) is

$$\boldsymbol{\theta}(\tau) = \boldsymbol{\theta}_0 + \int_{\tau_0}^{\tau} (-\boldsymbol{K}_\theta f_\theta) \mathrm{d}s , \qquad (28)$$

where $\boldsymbol{\theta}_0 = \boldsymbol{\theta}(\tau_0)$ is the initial condition. Now we introduce a limited integrator symbol as

$$y(t) = \overline{\underline{y(t_0) + \int_{t_0}^{t} x(s) \mathrm{d}s}}^{\,a}_{\,b} , \qquad (29)$$

where $a$ and $b$ denote the upper and the lower limit for the variable $y$, respectively. This symbol implies the differential relation as

$$\frac{\mathrm{d}y}{\mathrm{d}t} = \begin{cases} 0 & (y = a, x \geq 0) \\ 0 & (y = b, x \leq 0) \\ x & \text{others} \end{cases} . \qquad (30)$$

The limited integrator, acting on Eq. (15), will be employed to address the parameter bounds in Problem 1. Regarding its rationality, it is guaranteed by Theorem 2 in the following. First a Lemma is presented before the proof.

**Lemma 2:** The relation $\boldsymbol{K}_\theta^{1/2} \left( \boldsymbol{I} - (\tilde{\boldsymbol{h}}_\theta \boldsymbol{K}_\theta^{1/2})^+ (\tilde{\boldsymbol{h}}_\theta \boldsymbol{K}_\theta^{1/2}) \right) (\boldsymbol{K}_\theta^{1/2})^{-1} = \boldsymbol{I} - \tilde{\boldsymbol{h}}_\theta^+ \tilde{\boldsymbol{h}}_\theta$ holds when $\boldsymbol{K}_\theta$ is an $n \times n$ dimensional diagonal positive-definite matrix and $\tilde{\boldsymbol{h}}_\theta$ is the Jacobi matrix of $\tilde{\boldsymbol{h}}$ given in Theorem 1.

*Proof.* Due to the special form of $\tilde{\boldsymbol{h}}_\theta$ and $\boldsymbol{K}_\theta$, it is easy to verify that

$$(\tilde{\boldsymbol{h}}_\theta \boldsymbol{K}_\theta^{1/2})^+ = (\boldsymbol{K}_\theta^{1/2})^{-1} (\tilde{\boldsymbol{h}}_\theta)^+ . \qquad (31)$$

Thus, there is

$$\begin{aligned}
& \boldsymbol{K}_\theta^{1/2}\left(\boldsymbol{I} - (\tilde{\boldsymbol{h}}_\theta \boldsymbol{K}_\theta^{1/2})^+ (\tilde{\boldsymbol{h}}_\theta \boldsymbol{K}_\theta^{1/2})\right)(\boldsymbol{K}_\theta^{1/2})^{-1} \\
& = \boldsymbol{I} - \boldsymbol{K}_\theta^{1/2}(\tilde{\boldsymbol{h}}_\theta \boldsymbol{K}_\theta^{1/2})^+ (\tilde{\boldsymbol{h}}_\theta \boldsymbol{K}_\theta^{1/2})(\boldsymbol{K}_\theta^{1/2})^{-1} \\
& = \boldsymbol{I} - \tilde{\boldsymbol{h}}_\theta^+ \tilde{\boldsymbol{h}}_\theta
\end{aligned} \qquad (32)$$

Q.E.D. □

**Theorem 2**: *For Problem* 1, *solve the IVP defined by the following integral equation*

$$\boldsymbol{\theta} = \boldsymbol{\theta}_0 + \overbrace{\underbrace{\int_{\tau_0}^{\tau} (-\boldsymbol{K}_\theta f_\theta)\, \mathrm{d}t}_{\theta_l}}^{\theta_u}, \qquad (33)$$

where $\boldsymbol{\theta}_0 \in \mathbb{D}_f$ is arbitrary feasible initial condition, $\boldsymbol{\theta}_l$ is the Lower limit and $\boldsymbol{\theta}_u$ is the Upper limit. $\boldsymbol{K}_\theta$ is an $n \times n$ dimensional diagonal positive-definite gain matrix. When $\tau \to +\infty$, $\boldsymbol{\theta}$ will satisfy the optimality condition given by Eqs. (6) and (7).

*Proof.* Since problem 1 may be solved with the general scheme in the last section. Here we will prove that the DOE established by Eq. (33) is equivalent to Eq. (25), when the positive-definite matrix $\boldsymbol{K}_\theta$ is diagonal. The limited integrator acts on the differential equations of parameters. For $\theta_i$, its derivative is zero if the parameter bound constraint is active and the derivative satisfies

$$\begin{aligned}
(-\boldsymbol{K}_\theta f_\theta)_i \geq 0 & \quad \text{if } \theta_i = (\theta^u)_i \\
(-\boldsymbol{K}_\theta f_\theta)_i \leq 0 & \quad \text{if } \theta_i = (\theta^l)_i
\end{aligned}, \qquad (34)$$

where the subscript "$i$" denotes the $i$ th component of the vector. Any constraint on $\theta_i$ satisfying Eq. (34) belongs to the set of $\mathbb{I}_p$ defined by Eq. (20), or the derivative relation will be violated. Note that $\tilde{\boldsymbol{h}}_\theta^+ \tilde{\boldsymbol{h}}_\theta$ has a special form as demonstrated in Sec. 3.3. This then determines the equivalent differential form of Eq. (33) as

$$\frac{\mathrm{d}\boldsymbol{\theta}}{\mathrm{d}\tau} = -(\boldsymbol{I} - \tilde{\boldsymbol{h}}_\theta^+ \tilde{\boldsymbol{h}}_\theta) \boldsymbol{K}_\theta f_\theta, \qquad (35)$$

with $\boldsymbol{\theta}|_{\tau=0} = \boldsymbol{\theta}_0$. On the other hand, Eq. (25) may be re-presented as (See Ref. [36] for detail)

$$\frac{\mathrm{d}\boldsymbol{\theta}}{\mathrm{d}\tau} = -\boldsymbol{K}_\theta^{1/2}\left(\boldsymbol{I} - (\tilde{\boldsymbol{h}}_\theta \boldsymbol{K}_\theta^{1/2})^+ (\tilde{\boldsymbol{h}}_\theta \boldsymbol{K}_\theta^{1/2})\right) \boldsymbol{K}_\theta^{1/2} f_\theta. \qquad (36)$$

According to Lemma 2, we have that Eq. (35) and Eq. (36) are exactly the same. Therefore, when $\tau \to +\infty$, $\boldsymbol{\theta}$ will satisfy the optimality condition as the solution of Eq. (25) does. □

Upon the equivalence relation, the KKT multipliers may be estimated by Eq. (27), as

$$\boldsymbol{\pi} = -(\tilde{\boldsymbol{h}}_\theta \tilde{\boldsymbol{h}}_\theta^{\mathrm{T}})^{-1} \tilde{\boldsymbol{h}}_\theta f_\theta = -\tilde{\boldsymbol{h}}_\theta f_\theta. \qquad (37)$$

This means that the multiplier corresponding to the active bound constraint is $(\boldsymbol{\pi}_I)_i = -(f_\theta)_i$ for the upper bounds, and $(\boldsymbol{\pi}_I)_i = (f_\theta)_i$ for the lower bounds, which is consistent with the optimality conditions in Sec. 2. Besides the view of equivalence to Eq. (25). Actually, it may be shown that with the Eq. (35), the objective function always decreases for an $n \times n$ dimensional diagonal positive-definite gain matrix $\boldsymbol{K}_\theta$. See

$$\begin{aligned}\frac{\mathrm{d}f}{\mathrm{d}\tau} &= f_\theta^\mathrm{T} \frac{\mathrm{d}\theta}{\mathrm{d}\tau} = -f_\theta^\mathrm{T}(I - \tilde{h}_\theta^+ \tilde{h}_\theta) K_\theta f_\theta \\ &= -\left((I - \tilde{h}_\theta^+ \tilde{h}_\theta) f_\theta\right)^\mathrm{T} K_\theta f_\theta \leq 0\end{aligned} \quad (38)$$

Note that $K_\theta$ is restrict to be *diagonal* and cannot be an arbitrary positive-definite matrix. When $K_\theta$ is not diagonal, it may be occurred that

$$\frac{\mathrm{d}f}{\mathrm{d}\tau} = f_\theta^\mathrm{T} \frac{\mathrm{d}\theta}{\mathrm{d}\tau} = -f_\theta^\mathrm{T}(I - \tilde{h}_\theta^+ \tilde{h}_\theta) K_\theta f_\theta > 0. \quad (39)$$

For example, let $f_\theta = \begin{bmatrix} 1 \\ 2 \end{bmatrix}$, $K_\theta = \begin{bmatrix} 0.2 & -1 \\ -1 & 10 \end{bmatrix}$, $\tilde{h}_\theta = \begin{bmatrix} 0 & 1 \end{bmatrix}$. Then there is $\frac{\mathrm{d}f}{\mathrm{d}\tau} = 1.8$. From Eq. (35), it is shown that its form will not be altered by the strictness of the complementarity condition. The only difference is that the derivative takes the value of zero in a gradual or sudden way at the bounds. In particular, since the objective function $f$ is a valid Lyapunov function holding for arbitrarily large region around the minimums, Theorem 2 implies a global convergence property of the proposed approach.

**Remark 2:** Solve Problem 1 with the approach suggested in Theorem 2; the solution will converge to the optimal solution globally.

### 4.2. Numerical issues

In the numerical integration, it is impossible to detect the exact time when the bound constraint is active and there may be errors on the parameters when the limit integrator works. In order to eliminate the numerical error when applying the proposed unconstrained-like dynamic method, the differential form to the limited integrator, given by Eq. (30), is modified as

$$\frac{\mathrm{d}y}{\mathrm{d}t} = \begin{cases} -k^u (y-a) & (y \geq a, x \geq 0) \\ -k^l (y-b) & (y \leq b, x \leq 0) \\ x & \text{others} \end{cases}, \quad (40)$$

where $k^u$ and $k^l$ are the limiting gains. Naturally, when $y$ takes the value of the bound, Eq. (40) degenerates to Eq. (30).

As shown in Eq. (40), when forcing the derivative change from $x$ to follow the manually set first-order dynamics, two conditions needs to be satisfied. Loop-based computation may be time-consuming, especially for the large-scale problems, and special attention should be paid on the code optimization.

### 5. Examples

First, an illustrative example as follows is considered.

**Example 1:**

$$\begin{aligned} \min \quad & f(\boldsymbol{\theta}) = (\theta_1 + 1)^2 + (\theta_2 - 2)^2 \\ \text{subject to} \quad & 0 \leq \theta_1 \leq 10 \\ & 0 \leq \theta_2 \leq 10 \end{aligned}$$

For this problem, as shown in Fig.1, by investigating the contour of the objective function and the feasible solution region, it is easy to locate the optimal solution $\hat{\boldsymbol{\theta}} = \begin{bmatrix} \hat{\theta}_1 \\ \hat{\theta}_2 \end{bmatrix} = \begin{bmatrix} 0 \\ 2 \end{bmatrix}$.

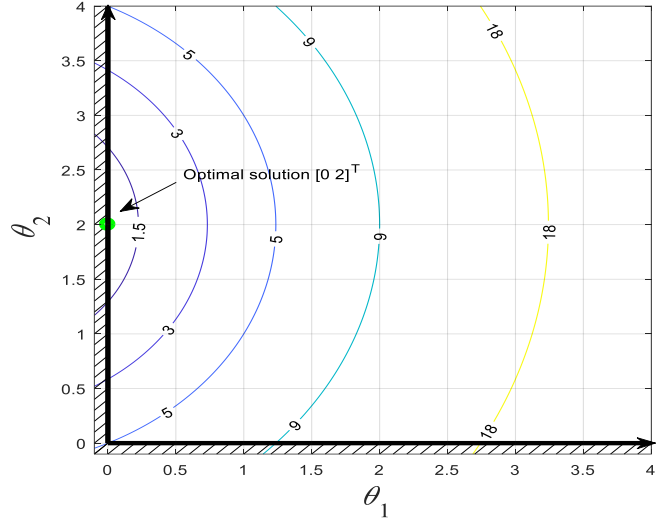

Figure 1. The optimal solution of Example 1 by geometric analysis.

Now we will use the proposed unconstrained-like method to seek the solution, and the limited integral equation is

$$\boldsymbol{\theta} = \boldsymbol{\theta}\big|_{\tau=0} + \overbrace{\int_{\tau_0}^{\tau} 2\boldsymbol{K}_{\boldsymbol{\theta}} \begin{bmatrix} \theta_1 + 1 \\ \theta_2 - 2 \end{bmatrix} \mathrm{d}s}^{[0\ 10]^{\mathrm{T}}}_{[0\ 10]^{\mathrm{T}}},$$

where the gain matrix is $\boldsymbol{K}_{\boldsymbol{\theta}} = \begin{bmatrix} 0.5 & 0 \\ 0 & 1 \end{bmatrix}$, and the initial guesses of the parameters are $\boldsymbol{\theta}\big|_{\tau=0} = \begin{bmatrix} 5 \\ 5 \end{bmatrix}$. Besides, the limiting gains in Eq. (40) are all set to be 1. To solve the resulting IVP, the ODE integrator 'ode45' in MATLAB, which implements explicit Runge-Kutta (4, 5) formula of the Dormand-Prince pair, was employed on an integration time horizon of 100s. The relative error tolerance was set as $1 \times 10^{-3}$ and the absolute error tolerance was set as $1 \times 10^{-6}$ in the integrator.

Figure 2 gives the solution of the parameters against the virtual time $\tau$. Both parameter solutions converge to the optimal solution as time advances. For comparison, the solution with the general dynamic method, as given in Sec. 3, is also presented. It is shown that the results are the same. However, when testing the computation time with both approaches, the time consumption is 0.014s for the unconstrained-like method and 0.049s for the general dynamic method, respectively. The proposed method is more efficient in that no intermediate QP sub-problem is involved. In Fig. 3, we gives the convergence trajectories in the parameter coordinate plane, obtained by the unconstrained-like method with different initial guesses. It is shown that the proposed approach succeeds to reach the optimal solution. When the parameter bound constraint becomes active, the motion trajectory changes the direction and continues to tend to the optimal solution.

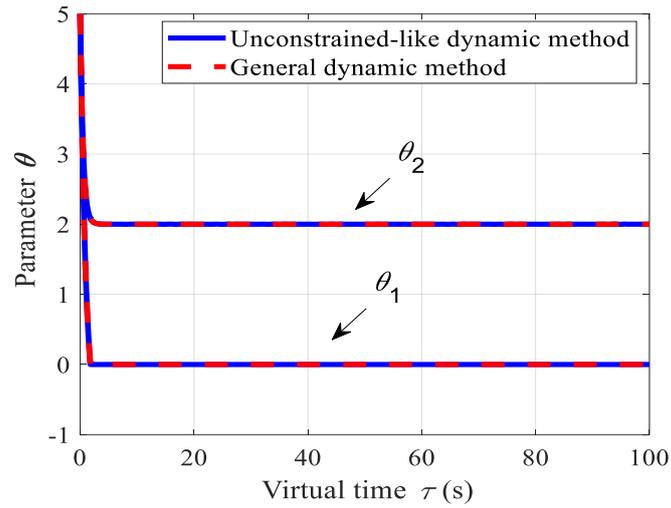

Figure 2. The dynamic motion curves of the parameters by the proposed unconstrained-like dynamic method and the general dynamic method.

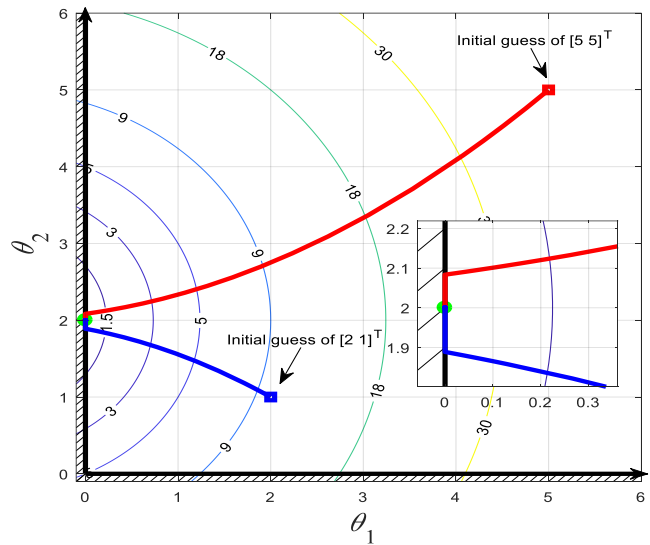

Figure 3. The convergence trajectories in the parameter coordinate plane from different initial conditions.

In particular, in order to examine the effect with a non-diagonal positive-definite gain matrix in the proposed method. The gain matrix is altered to be $\boldsymbol{K_\theta} = \begin{bmatrix} 0.5 & 0.2 \\ 0.2 & 1 \end{bmatrix}$ while the other setting is similar. Figure 4 gives the dynamic motion curves for the parameters. Now the right solution cannot be obtained and the parameters finally stay at the value of $\begin{bmatrix} 0 & 1.8 \end{bmatrix}^T$.

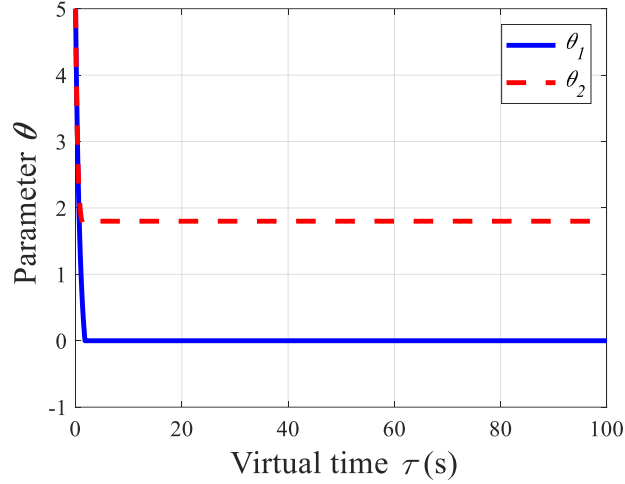

Figure 4. The dynamic motion curves of the parameters under non-diagonal positive-definite gain matrix

Next we consider a well-known problem where the strict complementarity condition does not hold [19]; see

**Example 2**:

$$\min\ f(\boldsymbol{\theta}) = -\frac{1}{2}(\theta_1^2 - \theta_2^2) - \theta_1^2 \theta_2 + \theta_1$$

subject to

$$0 \leq \theta_1 \leq 1$$
$$0 \leq \theta_2 \leq 1$$

For this problem, the minimizer point is $\hat{\boldsymbol{\theta}} = \begin{bmatrix} \hat{\theta}_1 \\ \hat{\theta}_2 \end{bmatrix} = \begin{bmatrix} 0 \\ 0 \end{bmatrix}$, and the gradient for $\hat{\theta}_2$ is zero. Then the strict complementarity is violated in this parameter. With the proposed unconstrained-like method, the limited integral equation is

$$\boldsymbol{\theta} = \boldsymbol{\theta}|_{\tau=0} + \int_{\tau_0}^{\tau} 2\boldsymbol{K}_{\boldsymbol{\theta}} \overline{\begin{bmatrix} -\theta_1 - 2\theta_1 \theta_2 + 1 \\ \theta_2 - \theta_1^2 \end{bmatrix}}_{[0\ 1]^T}^{[0\ 1]^T} \mathrm{d}s,$$

where the gain matrix is $\boldsymbol{K}_{\boldsymbol{\theta}} = \begin{bmatrix} 0.5 & 0 \\ 0 & 1 \end{bmatrix}$, and the initial condition is $\boldsymbol{\theta}|_{\tau=0} = \begin{bmatrix} 0.5 \\ 0.5 \end{bmatrix}$. The limiting gains are also set to be 1. Moreover, to compare the performance when the strict complementarity condition is satisfied, another set of bound constraints, which results in the optimal solution of $\hat{\boldsymbol{\theta}} = \begin{bmatrix} \hat{\theta}_1 \\ \hat{\theta}_2 \end{bmatrix} = \begin{bmatrix} 0.1 \\ 0.1 \end{bmatrix}$, is also considered; that is

$$0.1 \leq \theta_1 \leq 1$$
$$0.1 \leq \theta_2 \leq 1 \ .$$

In solving the resulting IVPs, the ODE integrator 'ode45' in MATLAB, with a relative error tolerance of $1\times 10^{-3}$ and an absolute error tolerance of $1\times 10^{-6}$, was employed on an integration time horizon of 50s.

The computation results are plotted in Figs. 5 and 6. In Fig. 5, it is shown that the solutions, for the box constraints of [0, 1] and [0.1, 1] both, converge rapidly to the optimal solution under the proposed method, no matter whether the complementarity condition are strictly satisfied or not. The computation time are 0.0179s and 0.0328s, respectively. Although the convergence of $\theta_2$ seems slower than $\theta_1$, it is interesting to find that the computation time, used for the case where the strict complementarity is violated, is shorter. This is because for the box constraint of [0, 1], $\theta_2$ approaches zero asymptotically, while for the box constraint of [0.1, 1], the slope has a sudden change and this is disadvantageous to the numerical integration. In Fig. 6, the motion trajectory on the parameter plane, where the contours for the objective function are demonstrated, is plotted. It shows how the optimal solution is achieved under the gradient dynamics.

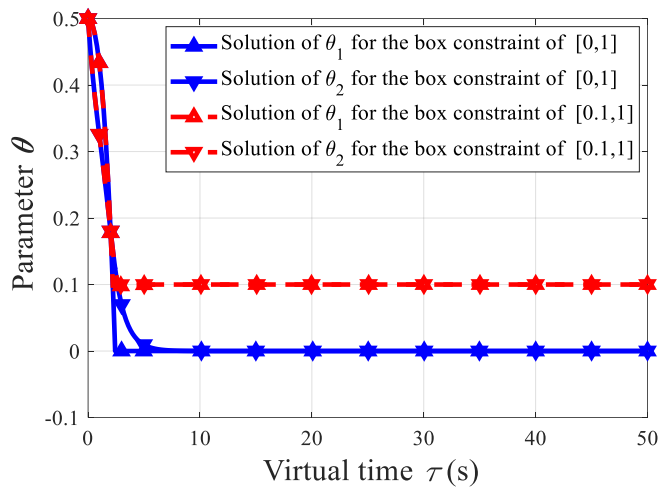

Figure 5. The dynamic motion curves of the parameters for the two box constraints.

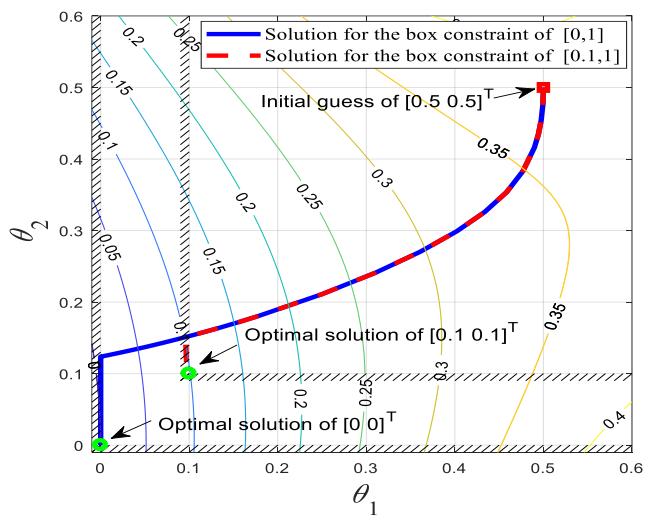

Figure 6. The convergence trajectories in the parameter coordinate plane for the two box constraints.

In the preceding examples, the number of the parameters is small, while it is interesting to investigate the efficiency of the proposed method for more parameters. Consider the classic Generalized Wood function (GENWOOD) optimization problem [14].

**Example 3:**

$$\min \ f(\boldsymbol{\theta}) = 1 + \sum_{i \in \Lambda} \begin{pmatrix} 100(\theta_{i+1} - \theta_i^2)^2 + (1-\theta_i)^2 + 90(\theta_{i+3} - \theta_{i+2}^2)^2 \\ +(1-\theta_{i+2})^2 + 10(\theta_{i+1} + \theta_{i+3} - 2)^2 + 0.1(\theta_{i+1} - \theta_{i+3})^2 \end{pmatrix}$$

subject to

$$1.1 \leq \theta_j \leq 2.1 \quad j = 1, 2, ..., n$$

where $\boldsymbol{\theta} \in \mathbb{R}^n$ is the optimization parameter vector, $n$ is a multiple of 4 and $\Lambda = \{1, 5, 9, ..., n-3\}$. The optimal solution for this example is $\hat{\boldsymbol{\theta}} = [1.1 \ 1.1753 \ 1.1 \ 1.1715 \ ... \ 1.1 \ 1.1753 \ 1.1 \ 1.1715]_{n \times 1}^T$. Through the proposed unconstrained-like method, the integral equation is given by Eq. (33), where $f_{\boldsymbol{\theta}} = \begin{bmatrix} \boldsymbol{d}_1^T & \boldsymbol{d}_2^T & ... & \boldsymbol{d}_j^T & ... & \boldsymbol{d}_{n/4}^T \end{bmatrix}_{n \times 1}^T$ and

$$\boldsymbol{d}_j = \begin{bmatrix} -400(\theta_{i+1} - \theta_i^2)\theta_i - 2(1-\theta_i) \\ 200(\theta_{i+1} - \theta_i^2) + 20(\theta_{i+1} + \theta_{i+3} - 2) + 0.2(\theta_{i+1} - \theta_{i+3}) \\ -360(\theta_{i+3} - \theta_{i+2}^2)\theta_{i+2} - 2(1-\theta_{i+2}) \\ 180(\theta_{i+3} - \theta_{i+2}^2) + 20(\theta_{i+1} + \theta_{i+3} - 2) - 0.2(\theta_{i+1} - \theta_{i+3}) \end{bmatrix}; j = 1, 2, ..., \frac{n}{4}; i = 4j - 3$$

$\boldsymbol{\theta}_l = [1.1 \ 1.1 \ ... \ 1.1]^T$ and $\boldsymbol{\theta}_u = [2.1 \ 2.1 \ ... \ 2.1]^T$. $\boldsymbol{K}_{\boldsymbol{\theta}}$ was set as $\boldsymbol{K}_{\boldsymbol{\theta}} = \boldsymbol{I}_{n \times n}$. Consider the case of $n = 100$. The initial conditions of $\boldsymbol{\theta}$ were given by $\boldsymbol{\theta}|_{\tau=0} = [1.1 \ 1.1 \ ... \ 1.1]_{100 \times 1}^T$. Limiting gains of 1 were applied. Because the dynamics of the DOE are fast and a little stiff, to enhance the numerical stability when solving the resulting IVPs, the stiff ODE integrator "ode15s" in MATLAB, which uses a variable-step, variable-order solver based on the numerical differentiation formulas (NDFs) of orders 1 to 5, was selected for the numerical integration along the time horizon of 10s. The relative error tolerance was $1 \times 10^{-3}$ and an absolute error tolerance was $1 \times 10^{-6}$.

Figure 7 gives the profiles of the parameters, showing that the solutions converge to the optimal very rapidly and the minimizer point is almost obtained after $\tau = 0.04$s. For comparison, other common iterative methods for the constrained problems were employed, including the active-set method, the trust-region-reflective method (one variant of trust-region method), the interior-point method and the Sequential Quadratic Programming (SQP) method, by invoking the function of "fmincon" in MATLAB. Besides, 50 computations starting from different initial conditions of $\boldsymbol{\theta}|_{\tau=0}$ were tried and their values were sampled from a uniform distribution within [1.1, 2.1]. All the methods succeed to get the optimal solution while their efficiency are different. The time consumption is presented in Table 1. It is shown that the proposed dynamic method seek the right solution with the highest efficiency. We also tried the unconstrained optimization method by parameter transformation, with the function "fminsearchbnd" [22]. However, it fails to get the solution.

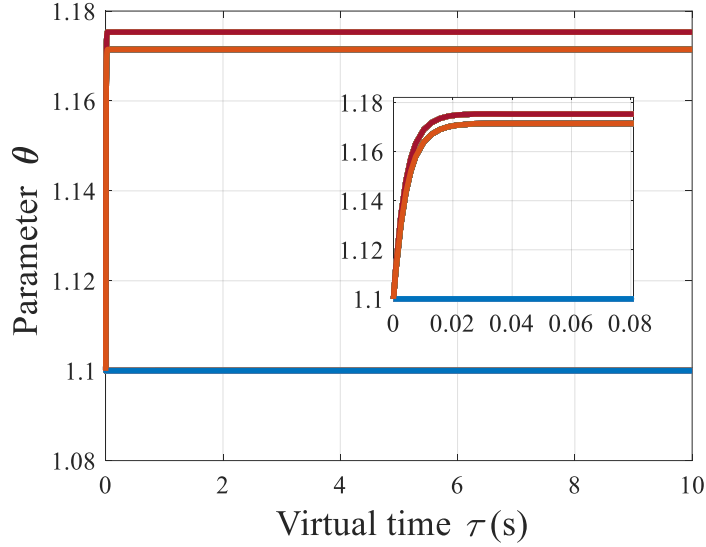

Figure 7. The dynamic motion curves of the parameters for example 3. There are 100 curves in total.

Table 1. Time consumption for Example 3 with different methods.

|  | The proposed dynamic method | The active-set method | The trust-region-reflective method | The interior-point method | The SQP method |
|---|---|---|---|---|---|
| Time consumption for the initial guess of [1.1 1.1 1.1 … 1.1]$^T$ (s) | 0.0746 | 1.5270 | 0.0811 | 0.5756 | 0.1284 |
| Average time consumption for the 50 computations (s) | 0.0893 | 1.6857 | 0.0878 | 0.6110 | 0.1403 |

Now a practical problem adapted from Ref. [36] is investigated, which is a Partial Differential Equation (PDE) constrained problem. The aim is to seek the thermal conductivity $k(T)$ for the heat shield material. Through parameterization, the problem may be formulated as

**Example 4**:

$$\min \ J(\boldsymbol{\theta}) = \int_{t_0}^{t_f} \sum_{i=1}^{m} \left\{ w_i \left[ T(x_i,t; \tilde{k}(T)) - \tilde{T}(x_i,t) \right]^2 \right\} dt$$

subject to

$$\rho C_p \frac{\partial T(x,t)}{\partial t} = \frac{\partial}{\partial x}\left( \tilde{k}(T) \frac{\partial T(x,t)}{\partial x} \right), \ \tilde{k}(T) = \sum_{i=1}^{n} (\sum_{j=1}^{i} \theta_j) \phi_i(T)$$

$$T(x,0) = T_{t0}, T(0,t) = T_{x0}, T(L,t) = T_{xL}$$

$$\tilde{T}(x_i,t) = T(x_i,t) + v(t) \ i = 1,2,...m$$

$$0 \le \theta_j \quad j = 1,2,...,n$$

where $T(x,t)$ is the temperature, $\rho$ is the density of material and $C_p$ is the specific heat. $\tilde{k}(T)$ is the modeling of $k(T)$. $\theta_j \ (j=1,2,...,n)$ are the optimization parameters. In particular, $\sum_{j=1}^{i} \theta_j \ (i=1,2,...,n)$ represents the thermal conductivity value at the temperature of $T_i$. The reason that we choose such modeling is to ensure the monotonic increase of the thermal conductivity. $\phi_i(T) = \prod_{j=1, j \neq i}^{n} \frac{(T-T_j)}{(T_i-T_j)}$ is the Lagrange interpolation polynomial. $\tilde{T}$ is the measured temperature and $v(t)$ is the measurement noise. $x_i$ is the position coordinate for the $i$-th measurement point and $m$ is the number of the measurement points. $w_i$ is the weight for the $i$-th measurement information and $\sum_{i=1}^{m} w_i = 1$. $[t_0, t_f]$ is the time span for measurement.

For the material studied, let $\rho$ = 1000 kg·m$^{-3}$, $C_p$ = 1000 J·kg$^{-1}$·K$^{-1}$, and $L$ = 10 mm. The initial condition for the heat conduction PDE was $T_{t0} = 600$ K. The boundary conditions were $T_{x0} = \begin{cases} (600+16t)\text{K} & t \leq 25s \\ 1000\text{K} & t > 25s \end{cases}$ and $T_{xL} = 600$ K. Two measurement points were available, and the noise obeyed a normal distribution with zero mean and a standard deviation of 6 K. The measurement times were $t_0$ = 0 s and $t_f$ = 100 s. To identify the thermal conductivity, $n = 5$ was used for the parameterization. The initial values of the optimization parameters were $\theta|_{\tau=0} = [2 \ 2 \ 2 \ 2 \ 2]^T$. Although the number of the optimization parameters is small, their motion scale is different and the resulting ODE is moderately stiff. Thus, the ODE integrator 'ode15s' is again employed for an integration time horizon of 200 s, with a relative error tolerance of $1 \times 10^{-3}$ and an absolute error tolerance of $1 \times 10^{-6}$.

In Fig. 8, the dynamic motion curves of the optimization parameters are presented. It is found that the parameters tend to the optimal solution rapidly and their values are almost fixed after $\tau$ = 100 s. During the computation, the bound constraint for $\theta_2$ and $\theta_3$ are active and they depart the bounds soon, as shown by the close-up. When $\tau$ = 200 s, we get $\theta|_{\tau=200} = [0.8614 \ 0.2129 \ 0.2176 \ 0.3563 \ 0.4465]^T$. In Fig. 9, the thermal conductivity identified, i.e., $\tilde{k}(T)$, are compared with the true quantity. They coincide well with each other, and the maximum error is only 0.0247 W·m$^{-1}$·K$^{-1}$.

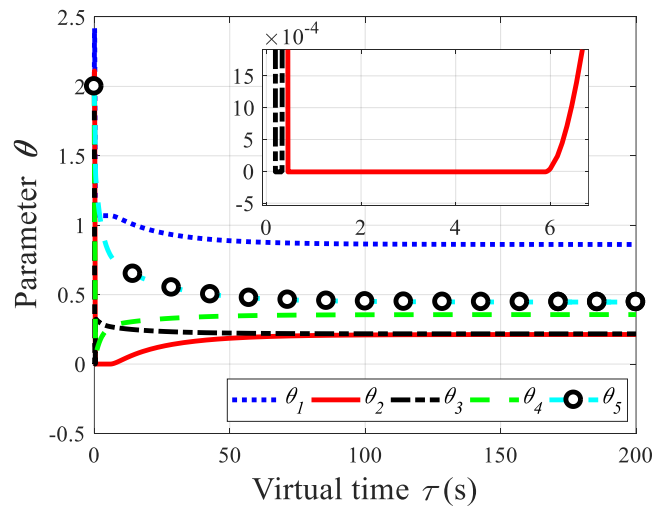

Figure 8. The dynamic motion curves of the discretized parameters for the thermal conductivity.

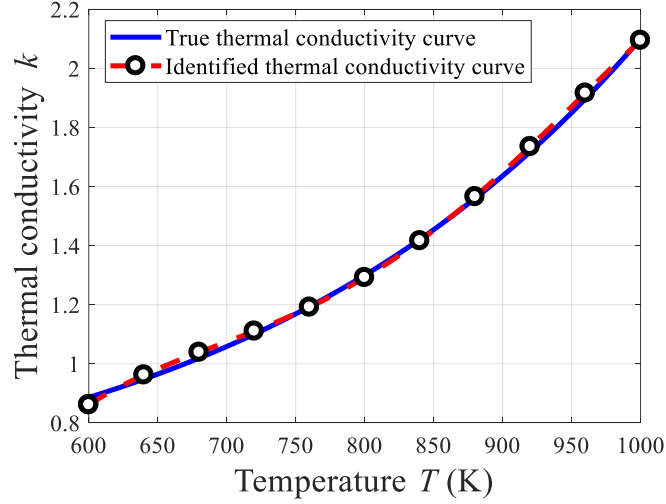

Figure 9. Comparison between the identified thermal conductivity and the true quantity.

To further test the proposed unconstrained-like dynamic method, 50 computations were carried out, with initial values sampled uniformly within the interval of [0, 2]. For comparison, this example is also computed with the general dynamic method and the common iterative methods, starting from the same initial values. Table 2 summarizes the results. In this example, failure occurs to the iterative methods when seeking the minimizer, while the proposed dynamic method has a superior performance overall. It not only produces the solutions fast, but also behaves robustly and all the computations find the optimal solution. The interior-point method costs the least time for this example. However, it encounters the failure in the computations, as the other iterative methods did. In addition, it is well-known that the scaling is important for the numerical optimization. In the numerical computation, it is found that the iterative methods are more likely to fail when using unreasonable scaling. However, for the dynamic method, it can always get the optimal solution upon the mature integration techniques.

Table 2. Computation results for Example 4 with different methods.

| | The unconstrained-like dynamic method | The general dynamic method | The active-set method | The trust-region-reflective method | The interior-point method | The SQP method |
|---|---|---|---|---|---|---|
| Average time consumption (s) | 60.20 | 60.96 | 114.64 | 71.34 | 21.35 | 107.69 |
| Failure times | 0 | 0 | 2 | 9 | 1 | 6 |

## 6. Conclusion

An easy-to-use unconstrained-like dynamic method for the problem with simple bounds is proposed, which relies on the dynamics for the unconstrained problem and simply employs the magnitude limiter to address the bounds. It is proved that the resulting Dynamic Optimization Equation (DOE) is equivalent to that from the general dynamic method when the gain matrix is diagonal, and thus its solution will meet the optimality condition with global convergence. Because the dynamics for the parameter that violates the strict complementarity condition is continuous, it is even more advantageous to the numerical integration. For the control-based view, even if only the first-order optimality conditions will be guaranteed, the solution will not halt at a maximum or a saddle (unless they are the initial guesses), since those solutions are not asymptotically stable equilibriums. Because the intermediate sub-problems used to address the inequality constraint during the computation process are not necessary anymore, the proposed method displays good performance, as shown by the numerical examples.

**References**

[1]   Gill, P. E., Murray, W., Wright, M. H.: *Practical Optimization*, Academic Press, New York (1981)

[2]   Nocedal, J., Wright, S. J.: *Numerical Optimization*, Springer Science, New York (2006)

[3]   Hager, W. W., Zhang, H.: Recent advances in bound constrained optimization. In F. Ceragioli, A. Dontchev, H. Furuta, K. Marti, & L. Pandolfi (Eds.), *System Modeling and Optimization*. Berlin, Europe: Springer. pp. 67–82 (2006)

[4]   Andrei, N.: Simple Bound Constraints Optimization. In: *Continuous Nonlinear Optimization for Engineering Applications in GAMS Technology*. Springer Optimization and Its Applications, vol. 121. Springer, Cham. (2017)

[5]   Bertsekas, D. P.: On the Goldstein-Levitin-Polyak gradient projection method. *IEEE Transactions on Automatic Control*, 21, 174–184 (1976)

[6]   Bertsekas, D. P.: Projected Newton methods for optimization problems with simple constraints. *SIAM Journal on Control and Optimization*, 20, 221–246 (1982)

[7]   Birgin, E. G., Mart ńez, J. M.: Large-scale active-set box-constrained optimization method with spectral projected gradients. *Comput Optim Appl* 23:101–125 (2002)

[8]   Hager, W. W., Zhang, H.: A new active set algorithm for box constrained optimization. *SIAM Journal on Optimization*, 17, 526–557 (2006)

[9]   Dost ál, Z.: Box constrained quadratic programming with proportioning and projections. *SIAM Journal on Optimization*, 7, 871–887 (1997)

[10] Lin, C. J., More, J. J.: Newton's method for large-scale bound constrained problems. *SIAM Journal on Optimization*, 9:1100–1127 (1999)

[11] Qi, L., Tong, X. J., Li, D. H.: Active-Set Projected Trust-Region Algorithm for Box-Constrained Nonsmooth Equations. *Journal of Optimization Theory and Applications*, vol. 120, 601–625 (2004)

[12] Cristofari, A., De Santis, M., Lucidi, S., Rinaldi, F.: A two-stage active-set algorithm for bound-constrained optimization. *Journal of Optimization Theory and Applications*, 172: 369–401 (2017)

[13] Conn, A. R., Gould, N. I. M., Toint, P. L.: Global convergence of a class of trust region algorithms for optimization with simple bounds. *SIAM Journal on Numerical Analysis*, 25(2), 433–460 (1988)

[14] Conn, A. R., Gould, N. I. M., Toint, P. L.: Testing a class of methods for solving minimization problems with simple bounds on the variables. *Mathematics of Computation*, 50, 399–430 (1988)

[15] Friedlander, A., Martinez, J. M., Santos, S. A.: A new trust region algorithm for bound constrained minimization, *Applied Mathematics and Optimization*, 30: 235-266 (1994)

[16] Lescrenier, M. L.: Convergence of trust region algorithms for optimization with bounds when strict complementarity does not hold. *SIAM Journal on Numerical Analysis*, 28: 476-495 (1991)

[17] Coleman, T.F., Li, Y.: On the convergence of interior–reflective Newton methods for nonlinear minimization subject to bounds. *Math. Program*. 67, 189–224  (1994)

[18] Coleman, T.F., Li, Y.: An interior trust region approach for nonlinear minimization subject to bounds. *SIAM J. Optim*. 6, 418–445 (1996)

[19] Heinkenschloss, M., Ulbrich, M., Ulbrich, S.: Superlinear and quadratic convergence of affine-scaling interior-point Newton methods for problems with simple bounds without strict complementarity assumption. *Math. Program.* 86, 615–635 (1999)

[20] Wei, W., Dai, H., Liang, W. T.: A novel projected gradient-like method for optimization problems with simple constraint. *Computational and Applied Mathematics*, 39, 168 (2020)

[21] Facchinei, F., Lucidi, S.: An unconstrained minimization approach to the solution of optimization problems with simple bounds. (1996)

[22] D'Errico, J.: fminsearchbnd, fminsearchcon (https://www.mathworks.com/matlabcentral/fileexchange/8277-fminsearchbnd-fminsearchcon), MATLAB Central File Exchange. Retrieved September 15 (2021)

[23] Courant, R.: Variational methods for the solution of problems of equilibrium and vibrations. *Bull. Amer. Math. Soc*. 49, 1-23 (1943)

[24] Botsaris, C. A.: Differential gradient methods. *J. Math. Anal. Appl.* 63, 177-198 (1978)

[25] Schaffler, S., Warsitz, H.: A trajectory-following method for unconstrained optimization. *J. Optim Theory Appl.* 67, 133-140 (1990)

[26] Mertikopoulos, P., Staudigl, M.: On the convergence of gradient-like flows with noisy gradient input. *SIAM J. Optim.* 28, 163-197 (2018)

[27] Snyman, J. A.: A new and dynamic method for unconstrained minimization. *Appl. Math. Model.* 6, 449-462 (1982)

[28] Alvarez, F., Perez, J. M.: A dynamical system associated with newton's method for parametric approximations of convex minimization problems. *Appl. Math. Optim.* 38, 193-217 (1998)

[29] Tanabe, K.: An algorithm for constrained maximization in nonlinear programming. *J. Oper. Res. Soc. Japan*. 17, 184-201 (1974)

[30] Yamashita, H.: A differential equation approach to nonlinear programming. *Math. Program.* 18, 155-168 (1980)

[31] Evtushenko, Y. G., Zhadan, V. G.: Stable barrier projection and barrier Newton methods in nonlinear programming. *Optim. Method. Softw.* 3, 237-256 (1994)

[32] Schropp, J.: A dynamical systems approach to constrained minimization. *Numer. Func. Anal. Opt.* 21, 537-551 (2000)

[33] Shikhman, V., Stein, O.: Constrained optimization: projected gradient flows. *J. Optim. Theory Appl.* 140, 117-130 (2009)

[34] Moguerza, J. M., Prieto, F. J.: Combining search directions using gradient flows. *Math. Program*. 96, 529-559 (2003)

[35] Ali, M. M., Oliphant, T. L.: A trajectory-based method for constrained nonlinear optimization problems. *J. Optim Theory Appl*, 177(2), 1-19 (2018)

[36] Zhang, S., Liao, F., Kong, Y. N., He K. F.: Transform the non-linear programming problem to the initial-value problem to solve. arXiv preprint, arXiv: 1709.02242 [cs.SY] (2017)

[37] Khalil, H.K.: *Nonlinear Systems*, Prentice Hall, New Jersey (2002)

[38] Shampine, L., Reichelt, M.: The MATLAB ODE Suite, *SIAM Journal on Scientific Computing*, 18, 1–22 (1997)